\documentclass[10pt,a4paper]{article}
\usepackage[a4paper,margin=2.5cm]{geometry}
\usepackage[utf8]{inputenc}
\usepackage[T1]{fontenc}
\usepackage{lmodern}
\usepackage{amsmath,amsfonts,amssymb,amsthm,mathtools}
\usepackage[parfill]{parskip}
\usepackage{microtype}
\usepackage{xcolor}
\newcommand{\Needspace}[1]{\par\begingroup
  \dimen0=\pagegoal\advance\dimen0 by -\pagetotal
  \ifdim\dimen0<#1\relax\newpage\fi\endgroup}
\usepackage[colorlinks=true,allcolors=blue]{hyperref}
\usepackage[capitalize,nameinlink]{cleveref}
\pdfoutput=1
\newtheorem{theorem}{Theorem}[section]
\newtheorem{prop}[theorem]{Proposition}
\newtheorem{lemma}[theorem]{Lemma}
\newtheorem{corollary}[theorem]{Corollary}
\theoremstyle{definition}

\theoremstyle{remark}

\crefname{prop}{Proposition}{Propositions}
\Crefname{prop}{Proposition}{Propositions}
\crefname{rem}{Remark}{Remarks}
\Crefname{rem}{Remark}{Remarks}
\numberwithin{equation}{section}

\newcommand{\Ff}{F}
\newcommand{\ZZ}{\mathbb{Z}}
\newcommand{\CC}{\mathbb{C}}
\newcommand{\rep}{\operatorname{Rep}}
\newcommand{\fm}{\mathfrak{m}}
\newcommand{\fl}{\mathfrak{l}}
\newcommand{\Lrep}{\mathrm{L}}
\newcommand{\St}{\operatorname{St}}
\newcommand{\bC}{C}
\newcommand{\bD}{D}

\newcommand{\cS}{\mathcal{S}}
\newcommand{\cZ}{\mathcal{Z}}
\newcommand{\cE}{\mathcal{E}}
\newcommand{\sm}{\mathrm{sm}}
\newcommand{\one}{\mathbf{1}}
\newcommand{\hra}{\hookrightarrow}
\newcommand{\sra}{\twoheadrightarrow}
\DeclareMathOperator{\Hom}{Hom}
\DeclareMathOperator{\Ext}{Ext}
\DeclareMathOperator{\GL}{GL}
\DeclareMathOperator{\rank}{rank}

\newcommand{\De}{\Delta}

\DeclareMathOperator{\soc}{soc}

\title{Higher theta lifts via BGG resolutions}
\author{Johannes Droschl}
\date{\today}

\begin{document}
\maketitle
\begin{abstract}
We compute higher theta lifts of ladder representations of general linear
groups over division algebras over non-archimedean local fields.
We achieve this by constructing BGG-type
coresolutions and proving that their terms are acyclic for the covariant
theta functor. We then compute the resulting complex explicitly,
determining the higher theta lifts and the degrees in which they
are nonzero.
\end{abstract}
\section{Introduction}
Let $F$ be a non-archimedean local field.
The theta correspondence, introduced by Howe in \cite{Howe1979}, associates
representations of one member of a reductive dual pair $(G,H)$ to
representations of the other. Its starting point is the Weil
representation $\omega$ \cite{Weil1964}. For an irreducible smooth
representation $\pi$ of $G$, the maximal $\pi$-isotypic quotient of
$\omega$ has the form $\pi\boxtimes\Theta(\pi)$. The representation
$\Theta(\pi)$ is called the big theta lift. Howe duality states that
its maximal semisimple quotient, called the small theta lift, is either
zero or irreducible, and that a nonzero small theta lift determines $\pi$.
We consider type II pairs, consisting of two general linear groups.
Howe duality in this setting is established in \cite{Minguez2008}.
The big theta lift need not be irreducible, and the theta functor need
not be exact. Following Adams, Prasad and Savin \cite{APS2017}, we
consider the higher lifts
\[
\Theta^r(\pi)
=\Ext_G^r(\omega,\pi)_{H\text{-}\sm}^{\vee}.
\]
For $r=0$ this recovers the big theta lift. 

We compute these higher theta lifts explicitly for ladder representations
on reduced-norm character lines of general linear groups over division
algebras. Their definition is
recalled in \cref{S:Seg}, following \cite{LM2014}.
We construct coresolutions of such ladder representations by
costandard modules whose higher theta lifts vanish.
These are inspired by the Bernstein--Gelfand--Gelfand (BGG) resolution
of a finite-dimensional irreducible representation of a complex
semisimple Lie algebra by Verma modules \cite{BGG1975}.
More precisely, if $\lambda$ is a dominant integral weight and $W$ is
the Weyl group, the term in degree $-r$ of the BGG resolution is
\[
 \bigoplus_{\ell(w)=r}\Sigma^{\mathcal O}(w\mathbin{\cdot}\lambda),
 \qquad w\mathbin{\cdot}\lambda=w(\lambda+\rho)-\rho.
\]
Here $\Sigma^{\mathcal O}(\mu)$ is the Verma module of highest weight $\mu$, and $\rho$
is the half sum of the positive roots. 
The cohomology of this complex
is concentrated in degree $0$, where it is the simple representation
of highest weight $\lambda$.
Geometrically, it can be constructed from the stratification of a flag
variety by Schubert cells.
For general linear groups over a non-archimedean
field, the Arakawa--Suzuki functor transports the BGG resolution to
ladder representations \cite[\S4]{BC2015}.
This transport was already established by Suzuki
\cite[Theorem~5.1.1]{Suzuki1998}, the corresponding Hecke-algebra
resolution in the ordinary BGG case goes back to Cherednik
\cite{Cherednik1987}, with related work of Zelevinsky
\cite{Zelevinsky1987}, as recalled in
\cite[Remark~5.1.2]{Suzuki1998}.
Its Euler characteristic recovers Tadi{\'c}'s determinantal formula for
Speh representations \cite{Tadic1995}, extended to ladder representations
by Lapid and M\'inguez \cite[Theorem~1]{LM2014}.
In \cite[Theorem~1.2(5)]{RegularLM2018},
the authors extend this formula to representations attached to \emph{balanced}
multisegments.
This construction has a geometric origin. For a balanced multisegment,
the corresponding Schubert variety is smooth along the Schubert cell
determined by endpoint admissibility
\cite[Theorem~1.2(4)]{RegularLM2018}. Equivalently, the
Kazhdan--Lusztig polynomials on the admissible Bruhat interval are all
equal to $1$.
The Arakawa--Suzuki functor, followed by the corresponding Hecke-category
equivalences, transports the resulting complex to representations of
$p$-adic general linear groups.

After constructing a costandard coresolution of the representation in
question, we calculate its higher theta lifts as follows.
We first prove that the terms are acyclic for $\mathcal E^\bullet$, and then compute their theta lifts and the maps induced by the
differentials. Finally, we compute the resulting cohomology explicitly.
In type II the zeta integrals make the calculation of these maps
particularly explicit.
Recall that in that case, the Weil representation is realized on the space of
matrices. In equal rank this is also the space on which the
Godement--Jacquet zeta integrals are defined \cite{GJ1972}.
Dividing these integrals by the local $L$-factor gives maps which can
be evaluated at $s=1/2$ in the normalization used below. Their relation to the big and higher theta
lifts is studied in \cite{FSX2018,CLLTZ2025}. In our computation,
the ratio of the $L$-factors of two costandard modules determines
whether a differential survives after applying theta.
We deal first with representations in equal rank and then treat larger targets, using the
Kudla filtration and derivatives.

Let us now describe the results more explicitly.

Let $D$ be a central division algebra over $F$, of degree $e$, and
put $G_n=\GL_n(D)$. We let
$\nu_D=|\operatorname{Nrd}|^e$. Let $\chi$ be an unramified character
of $F^\times$, pulled back to each $G_n$ by the reduced norm.
A segment $[b,a]_\chi$, with $b-a\in\ZZ_{\geq0}$, has descending
cuspidal sequence
$\chi\nu_D^b,\ldots,\chi\nu_D^a$. A multisegment of the form
\[
 \fl=[b_1,a_1]_\chi+\cdots+[b_k,a_k]_\chi,\qquad
 a_1>\cdots>a_k,\quad b_1>\cdots>b_k,
\]
with all endpoint differences integral is called a \emph{ladder}.
Its total length is
$n=\sum_{r=1}^k(b_r-a_r+1)$. The Langlands classification attaches
to it an irreducible representation $\Lrep(\fl)$ of $G_n$.

In equal rank, the most interesting case occurs precisely when the lower
endpoint $\nu_D^{1/2}$ and the upper endpoint $\nu_D^{-1/2}$ both
occur in $\fl$. Write their indices as $i$ and $j$, respectively, and set
$d=j-i$. Deleting these two endpoints and pairing the remaining
endpoints in order gives a ladder $\fl^*$ of the same length.
By \cref{thm:main}, for $r>0$,
\[
 \Theta_n^r(\Lrep(\fl))\simeq
 \begin{cases}
  \Lrep(\fl^*)^\vee&r=d\text{ or }r=d-1,\\
  0&\text{otherwise}.
 \end{cases}
\]
If the two endpoints do not both occur, all positive lifts vanish.
The ordinary lift is $\Lrep(\fl)^\vee$, except when $d=1$. In that
case it is the nonsplit extension in \eqref{eq:Sigma}.
Thus the positive higher lifts occupy at most two consecutive degrees
and, when both occur, are isomorphic.

The same method applies in equal rank to the representations
attached to balanced multisegments of Lapid--M\'inguez
\cite[\S6]{RegularLM2018}, see \cref{gl:regular:theta}.

\paragraph{Acknowledgements.}
The author wants to thank Alberto M\'inguez and Corina Ciobotaru for the discussions surrounding this paper. This research was supported by a research grant (VIL53023) from VILLUM FONDEN.

\paragraph{Disclosure of AI assistance.}
The author supplied the idea of using BGG resolutions to compute
higher theta lifts in type II, however, all technical steps were developed in close cooperation with ChatGPT. Moreover, AI was used throughout the process of writing this manuscript, including drafting, adversarial proof checking, and language editing.
The author is responsible for the final version of the paper.

\section{General preliminaries}\label{sec:general-preliminaries}

\subsection{Representations, parabolic induction, and Jacquet modules}
\label{prelim:conventions}

Let $F$ be a non-archimedean local field
with absolute value $|-|$ and residue cardinality $q$. Fix a
uniformizer $\varpi$, so $|\varpi|=q^{-1}$.
For a reductive algebraic group $\mathbf G$ over $F$, put
$G=\mathbf G(F)$. Representations are smooth and have complex
coefficients. We write $\rep(G)$ for their category and compute
$\Ext_G$ in this category. The notation $\pi^\vee$ denotes the smooth
contragredient, which is exact on admissible representations of finite
length.
Fix a central division algebra $D$ over $F$ and put $e^2=[D:F]$.
For general linear groups over $D$ we use
$\nu_D(g)=|\operatorname{Nrd}(g)|^e$.
We will often just write $\nu$, omitting $D$.
For a parabolic subgroup $P=MU$ of one of the reductive groups used
below, with opposite $\overline P$ and modulus character $\delta_P$,
normalized parabolic induction is as usual
\[
 i_P^G(\sigma)
 =\{f:G\longrightarrow V_\sigma:\ f\text{ is smooth and }
 f(mug)=\delta_P(m)^{1/2}\sigma(m)f(g)\},
\]
with $G$ acting by right translation. The normalized Jacquet functor
is $r_P(\pi)=\delta_P^{-1/2}\otimes\pi_U$. These functors are exact,
and their adjunctions extend to Ext:
\begin{align*}
 \Ext_G^j(\pi,i_P^G\sigma)&\simeq\Ext_M^j(r_P\pi,\sigma),\\
 \Ext_G^j(i_P^G\sigma,\pi)&\simeq
       \Ext_M^j(\sigma,r_{\overline P}\pi).
\end{align*}
For these standard facts we refer to \cite{BZ1977} and Bernstein's
second adjointness theorem \cite{Bernstein1987}, the Ext formulations
are recorded in \cite[Lemma~3.2]{APS2017}.

\subsection{Descending segments}\label{prelim:segments}

Let $A=F$ or $D$, and let $\rho$ be an irreducible supercuspidal
representation of $\GL_h(A)$. Write
$\nu_\rho=|\operatorname{Nrd}|^{s_\rho}$, where $s_\rho>0$ is the
unique number for which $\rho\nu_\rho\times\rho$ is reducible
\cite[\S3.1]{RegularLM2018}. For $A=F$ one has $s_\rho=1$.
For $\rho=\chi\circ\operatorname{Nrd}$ on $D^\times=G_1$, where
$\chi$ is a character of $F^\times$, one has
$s_\rho=e$ and hence $\nu_\rho=\nu_D$.
For $u,v\in\mathbb{R}$ with
$u-v\in\ZZ_{\geq0}$, write
\[
 [u,v]_\rho=(\rho\nu_\rho^u,\rho\nu_\rho^{u-1},\ldots,\rho\nu_\rho^v),\qquad
 \De_\rho[u,v]=
   \soc(\rho\nu_\rho^u\times\rho\nu_\rho^{u-1}\times\cdots\times\rho\nu_\rho^v).
\]
The first is a formal segment, with upper endpoint $\rho\nu_\rho^u$,
lower endpoint $\rho\nu_\rho^v$ and length $u-v+1$. The second is its
essentially square-integrable representation on $\GL_{h(u-v+1)}(A)$
\cite[\S9]{Zelevinsky1980}, \cite[\S3.1]{RegularLM2018}.
When $\rho$ is unitary,
$\De_\rho[u,v]\nu_\rho^{-(u+v)/2}$ is unitary.
We recall also
\[
 [u,v]_\rho^\vee=[-v,-u]_{\rho^\vee},\qquad
 \De_\rho[u,v]^\vee\simeq\De_{\rho^\vee}[-v,-u].
\]
The formal segment $[u,u+1]_\rho$ is empty and we define
$\De_\rho[u,u+1]=\one_0$.

\subsection{Covariant theta and acyclic coresolutions}

For each specified dual pair $(G,H)$ and Weil representation $\omega$,
we compute the covariant functors
\[
 \mathcal E^j(\Pi)=\Ext_G^j(\omega,\Pi)_{H\text{-}\sm}
\]
and take the smooth contragredient afterwards to obtain
$\Theta^j(\Pi)=\mathcal E^j(\Pi)^\vee$. The subscript denotes smooth
vectors for $H$. The functor $\mathcal E^0$ is covariant and left
exact, since $\Pi$ occupies the second variable of $\Hom_G(\omega,\Pi)$.
The functors $\mathcal E^j$ form a cohomological $\delta$-functor
\cite[Lemma~7.4]{APS2017}. Taking the smooth contragredient reverses
arrows. We therefore use a coresolution to compute
$\mathcal E^\bullet$ and dualize its cohomology afterwards.

A resolution and a coresolution of $\Pi$ are respectively exact complexes
\[
 \cdots\longrightarrow C_1\longrightarrow C_0\longrightarrow \Pi
 \longrightarrow0,\qquad
 0\longrightarrow \Pi\longrightarrow D^0\longrightarrow D^1
 \longrightarrow\cdots .
\]
The coresolution is a cochain complex in nonnegative degrees.
For a covariant cohomological $\delta$-functor $\mathcal F^\bullet$,
an object $\Sigma$ is \emph{acyclic} if $\mathcal F^j(\Sigma)=0$ for $j>0$.

\begin{lemma}\label[lemma]{lem:acyclic-coresolution}
If $0\to \Pi\to D^0\to\cdots\to D^N\to0$ is a bounded
coresolution by acyclic objects, then naturally
\begin{equation}\label{eq:acyclic-coresolution}
 \mathcal F^j(\Pi)\simeq H^j\bigl(\mathcal F^0(D^\bullet)\bigr).
\end{equation}
The differentials on the right are induced by the given maps.
\end{lemma}
This follows by dimension shifting along the successive cokernels
of the coresolution.

\subsection{Category \texorpdfstring{$\mathcal O$}{O} and perverse sheaves}
\label{prelim:category-o}

Let $\mathfrak g=\mathfrak{gl}_k(\CC)$, with diagonal Cartan
$\mathfrak h$ and upper triangular Borel $\mathfrak b$.
Category $\mathcal O$ consists of finitely generated
$U(\mathfrak g)$-modules which are semisimple over $\mathfrak h$
and locally finite over $\mathfrak b$. For a weight $\mu$, put
\[
 \Sigma^{\mathcal O}(\mu)
 =U(\mathfrak g)\otimes_{U(\mathfrak b)}\CC_\mu.
\]
This Verma module has a unique simple quotient
$\pi^{\mathcal O}(\mu)$. Its costandard module
$\Upsilon^{\mathcal O}(\mu)$ is the direct sum of the duals of its
finite-dimensional weight spaces, with action twisted by the transpose
anti-involution. Write $\mathcal O_{\mu_-}$ for the block with the
generalized infinitesimal character of a regular integral antidominant
weight $\mu_-$. Subscripts $\tau\in S_k$ indicate highest weight
$\tau\mathbin{\cdot}\mu_-$.

Let $\mathcal B=\GL_k(\CC)/B$, with Schubert cells
$C_w\simeq\mathbb A^{\ell(w)}$, closures $X_w$ and inclusions $j_w$.
Write $D^b_{\mathrm{Sch}}(\mathcal B)$ for the ordinary bounded derived
category of complexes constructible along these cells, and
$\operatorname{Perv}_{\mathrm{Sch}}(\mathcal B)$ for its perverse heart.
Let $\mathsf T:\mathcal O_{\mu_-}\xrightarrow{\sim}\mathcal O_0$
be the translation equivalence with the principal block
\cite[\S7.1 and the theorem in \S7.8]{Humphreys2008}.
Beilinson--Bernstein localization \cite{BB1981}, together with the
work of Soergel \cite{Soergel1986,Soergel1989}, gives the geometric
realization formulated in \cite[Proposition~3.5.2]{BGS1996} as an equivalence
\[
 \mathsf L_B:D^b(\mathcal O_0)
       \xrightarrow{\ \sim\ }D^b_{\mathrm{Sch}}(\mathcal B)
\]
which identifies the usual $t$-structure with the perverse
$t$-structure. Thus the composite
\[
 \mathsf G_{\mu_-}=\mathsf L_B\circ D^b(\mathsf T):
 D^b(\mathcal O_{\mu_-})
       \xrightarrow{\ \sim\ }D^b_{\mathrm{Sch}}(\mathcal B)
\]
restricts to an exact equivalence
$\mathcal O_{\mu_-}\simeq\operatorname{Perv}_{\mathrm{Sch}}(\mathcal B)$.
Its values on the standard, costandard and simple objects are
\cite[Propositions~3.5.2 and~3.5.7]{BGS1996}
\[
 \begin{aligned}
 \mathsf G_{\mu_-}(\Sigma_\tau^{\mathcal O})
   &\simeq
      j_{\tau^{-1}!}\CC_{C_{\tau^{-1}}}[\ell(\tau)],\\
 \mathsf G_{\mu_-}(\Upsilon_\tau^{\mathcal O})
   &\simeq
      Rj_{\tau^{-1}*}\CC_{C_{\tau^{-1}}}[\ell(\tau)],\\
 \mathsf G_{\mu_-}(\pi_\tau^{\mathcal O})
   &\simeq \operatorname{IC}(X_{\tau^{-1}}).
 \end{aligned}
\]
The inverse appears because \cite[\S1.1]{BGS1996} labels the simple
attached to $X_w$ by $w^{-1}w_0\mathbin{\cdot}0$.
Here $\operatorname{IC}(X_w)$ is the intermediate extension of
$\CC_{C_w}[\ell(w)]$. For $\Pi,\Upsilon\in\mathcal O_{\mu_-}$,
the derived equivalence gives
\[
 \Ext^r_{\mathcal O_{\mu_-}}(\Pi,\Upsilon)
 \simeq\Hom_{D^b_{\mathrm{Sch}}(\mathcal B)}
 \bigl(\mathsf G_{\mu_-}(\Pi),\mathsf G_{\mu_-}(\Upsilon)[r]\bigr).
\]
This is the identification used below to compute Ext groups by
sheaf adjunction.

\subsection{The Arakawa--Suzuki functor}\label{prelim:arakawa-suzuki}

Put $V=(\CC^k)^{\otimes N}$ and let $\mathfrak n^-$ be the strictly
lower triangular Lie algebra. The Arakawa--Suzuki functor
\cite{AS1998,Suzuki1998} is
\[
 F_{\lambda,V}(\Pi)
 =H_0(\mathfrak n^-,\Pi\otimes V)_\lambda
 =\bigl((\Pi\otimes V)/\mathfrak n^-(\Pi\otimes V)\bigr)_\lambda.
\]
Thus one takes coinvariants for $\mathfrak n^-$ and then the
$\lambda$-weight space. Permutations of the tensor factors and Casimir
operators give this space an action of the graded affine Hecke algebra
$\mathbb H_N$ of type $A_{N-1}$
\cite[equation~(4.1.2) and Lemma~4.2.2]{BC2015}.
The functor is exact when $\lambda+\rho_k$ is dominant, where
$\rho_k$ is the half sum of positive roots \cite{AS1998}, see also
\cite[Remark~4.1.1]{BC2015}. Verma modules map to standard Hecke
modules or zero. The simple images and their nonvanishing conditions
are given in \cite[Theorems~3.2.1--3.2.2]{Suzuki1998}, we use the
sign convention and formulation of \cite[Theorem~4.2.3]{BC2015}.

On the central characters used below, the reduction from affine to
graded Hecke algebras and the Iwahori equivalence identify these Hecke
modules with representations of a split general linear group
\cite[\S1 and \S4]{BC2015}. The exact cuspidal-line equivalence of
\cite[Theorem~3.8]{RegularLM2018} then gives representations of
$\GL_N(D)$ on a fixed line generated by $\chi\circ\operatorname{Nrd}$,
where $\chi$ is a smooth character of $F^\times$. After absorbing a common
endpoint shift into $\chi$, write $\rep_N(\chi)$ for the full subcategory
of finite-length representations whose cuspidal support lies in
$\{\chi\nu_D^j:j\in\ZZ\}$. Thus the passage is
\[
 \operatorname{Perv}_{\mathrm{Sch}}(\mathcal B)
 \xleftarrow{\ \sim\ }\mathcal O_{\mu_-}
 \xrightarrow{\ F_{\lambda,V}\ }\mathbb H_N\text{-}\mathrm{mod}
 \longrightarrow\rep_N(\chi),
\]
with the last arrow understood on the relevant central characters.
The first equivalence computes the Ext coefficients used below.
The exact functor to $p$-adic representations transports the resulting
complex, including its differentials. It need not preserve those Ext
groups or give an equivalence on all of $\mathcal O_{\mu_-}$.

\section{General linear notation and segments}
\subsection{General linear induction}

Put $G_n=\GL_n(D)$ and let $\rep_n$ denote its category of smooth,
finite-length complex representations. We write
$\nu_n(g)=|\operatorname{Nrd}(g)|^e$, write $\nu=\nu_D$, and omit
the rank subscript when the group is clear.
A character $\chi$ of $F^\times$ is understood to be pulled back by
the reduced norm. Put $G_0=\{1\}$, with trivial representation $\one_0$.

For $n=n_1+\cdots+n_t$, let $P_{n_1,\ldots,n_t}=MU$ be the upper block
triangular parabolic subgroup of $G_n$, with Levi subgroup
$M\simeq G_{n_1}\times\cdots\times G_{n_t}$. We use the notation
\[
 \pi_1\times\cdots\times\pi_t
 =i_{P_{n_1,\ldots,n_t}}^{G_n}
   (\pi_1\boxtimes\cdots\boxtimes\pi_t)
\]
for the normalized induction fixed in the general preliminaries.

\subsection{Segments and standard modules}\label{S:Seg}

A \emph{multisegment} is a finite formal sum of segments. Length,
cuspidal support and duality extend to such sums.
Two formal segments $\mathfrak s_1,\mathfrak s_2$ are \emph{linked}
if their union is a segment and neither contains the other.
For a multisegment
$\fm=\sum_{i=1}^k[b_i,a_i]_{\rho_i}$, put
$\delta_i=\De_{\rho_i}[b_i,a_i]$ and let $n_i$ be its group rank.
The real central exponent $c_i$ is characterized by
$\delta_i\nu_{n_i}^{-c_i}$ being unitary.
Number the factors so that $c_1\geq\cdots\geq c_k$. The standard
and costandard modules are respectively
\[
 \St(\fm)=\delta_1\times\cdots\times\delta_k,\qquad
 \St'(\fm)=\delta_k\times\cdots\times\delta_1.
\]
The former has a unique irreducible quotient and the latter a unique
irreducible subrepresentation, both denoted by $\Lrep(\fm)$.
This is the Langlands classification in multisegment notation, due to
Bernstein--Zelevinsky and Zelevinsky in the split case
\cite{BZ1977,Zelevinsky1980}, and to Tadi{\'c} for inner forms in
characteristic zero \cite{Tadic1990}. A uniform treatment valid in
arbitrary characteristic is given in \cite[Appendix~A]{LM2016}.
For the conventions used here, see \cite[\S4.3]{EmertonHelm2014}
and the transfer in \cite[Theorem~3.8]{RegularLM2018}. In particular,
$\Lrep(\fm)^\vee\simeq\Lrep(\fm^\vee)$.

A multisegment $\fl=[b_1,a_1]_\chi+\cdots+[b_k,a_k]_\chi$ is called a
\emph{ladder} if all endpoint differences are integral,
$a_1>\cdots>a_k$ and $b_1>\cdots>b_k$ \cite{LM2014}.
Given $w\in S_k$, we write
\[
 \fl_w=\sum_{r=1}^k[b_{w(r)},a_r]_\chi,\qquad
 W(\fl)=\{w\in S_k:a_r\leq b_{w(r)}+1\text{ for all }r\}.
\]
Here $[a-1,a]_\chi$ is an empty segment, contributing $\one_0$ to an
induced product. A pair $[b,a]_\chi$ with $a>b+1$ makes the entire term
zero, as in \cite[Theorem~1]{LM2014}. Put
\begin{align*}
 \Sigma_w(\fl)&=\De_\chi[b_{w(1)},a_1]\times\cdots\times
              \De_\chi[b_{w(k)},a_k],\\
 \Upsilon_w(\fl)&=\De_\chi[b_{w(k)},a_k]\times\cdots\times
              \De_\chi[b_{w(1)},a_1].
\end{align*}
Thus $\Sigma_w(\fl)=\Upsilon_w(\fl)=0$ for $w\notin W(\fl)$. Although $\fl_w$ need
not be a ladder, its nonzero terms are standard and costandard modules.
\section{Acyclicity of costandard modules}
\subsection{The Weil representation and its rank filtration}

Let $M_{n,m}=M_{n,m}(D)$ be the space of $n\times m$ matrices over $D$. We let
$\omega_{n,m}$ be the space $\cS(M_{n,m})$ of locally constant, compactly
supported functions, with the smooth $G_n\times G_m$-action
\begin{equation}\label{eq:weil}
 (\omega_{n,m}(g,h)\varphi)(X)
 =\nu_n(g)^{-m/2}\nu_m(h)^{n/2}\varphi(g^{-1}Xh).
\end{equation}
We recall its rank filtration, see for example \cite[\S\S2--3]{Minguez2008}. Set $a=\min(n,m)$ and
\[
 U_r=\{X:\rank X\geq r\},\qquad \Omega_r=\{X:\rank X=r\}.
\]
Extension by zero identifies $\Gamma_r=\cS(U_r)$ with a subrepresentation of
$\omega_{n,m}$, and restriction to $\Omega_r$ gives
\[
 0=\Gamma_{a+1}\subset \Gamma_a\subset\cdots\subset \Gamma_0=\omega_{n,m},
 \qquad \Xi_r=\Gamma_r/\Gamma_{r+1}\simeq\cS(\Omega_r).
\]
Let $P_r=P_{r,n-r}$, $Q_r=P_{r,m-r}$, and let $\Pi_r^{\mathrm{reg}}$ be the regular
representation on $\cS(G_r)$, with action
$(\Pi_r^{\mathrm{reg}}(g,h)f)(x)=f(g^{-1}xh)$. Then
\begin{equation}\label{eq:rank-piece}
 \Xi_r\simeq i_{P_r\times\overline Q_r}^{G_n\times G_m}
              (\Pi_r^{\mathrm{reg}}\otimes\xi_r),
\end{equation}
where $\Pi_r^{\mathrm{reg}}$ is trivial on the remaining Levi factors, and
\begin{align*}
 \xi_r(g_r,g_{n-r};h_r,h_{m-r})
 ={}&\nu_r(g_r)^{(r-n-m)/2}\nu_{n-r}(g_{n-r})^{(r-m)/2}\\
 &\cdot\nu_r(h_r)^{(n+m-r)/2}\nu_{m-r}(h_{m-r})^{(n-r)/2}.
\end{align*}
When $n\leq m$, the full-rank piece gives
\begin{equation}\label{eq:open-piece}
 \Ext^j_{G_n}(\Xi_n,\pi)_\sm\simeq
 \begin{cases}\one_{m-n}\times\pi&j=0,\\0&j>0,\end{cases}
\end{equation}
by the regular-bimodule calculation of
\cite[Lemma~2.4 and equation~(14)]{APS2017}. Here $\one_h$ denotes the trivial
representation of $G_h$.

\subsection{Higher theta lifts of costandard modules}

Following \cite{APS2017}, we define the higher theta lifts of
$\pi\in\rep_n$ by
\begin{equation}\label{eq:theta-definition}
 \Theta^j_m(\pi)=\Ext^j_{G_n}(\omega_{n,m},\pi)_\sm^\vee,
 \qquad \cE^j_{n,m}(\pi)=\Ext^j_{G_n}(\omega_{n,m},\pi)_\sm.
\end{equation}
We also write $\Theta_m=\Theta_m^0$.
In equal rank we write $\cE^j_n=\cE^j_{n,n}$.

\begin{lemma}\label[lemma]{gl:rank-homology}
For $n\leq m$ and finite-length $\pi$, the representations
$\cE^j_{n,m}(\pi)$ have finite length and vanish for $j>n$. Moreover,
\[
 \sum_j(-1)^j[\cE^j_{n,m}(\pi)]=[\one_{m-n}\times\pi].
\]
\end{lemma}
\begin{proof}
For $D=F$, this is
\cite[Theorem~5.13 and Proposition~5.17]{APS2017} and the same argument works for a division algebra.
\end{proof}
We first record an induction formula valid for every $m=n+t$ with
$t\geq0$. Let $\delta=\De_\rho[u,v]$ be a segment representation
of $G_h$, with $h\leq n$, and let $\pi_0\in\rep_{n-h}$.
If the lower endpoint of $\delta$ is not $\nu^{(t+1)/2}$, then
for every $j\geq0$,
\begin{equation}\label{gl:segment-pullout}
 \cE^j_{n,m}(\delta\times\pi_0)
 \simeq\delta\times\cE^j_{n-h,m-h}(\pi_0).
\end{equation}
For $D=F$, this is \cite[Corollary~3.5]{CLLTZ2025}.
Over $D$, apply the same argument using the Kudla filtration of
Jacquet modules \cite{Kudla1986}, in the type-II form of
\cite[Proposition~3.2]{Minguez2008}.

We now specialize to equal rank.
\begin{lemma}\label[lemma]{lem:acyclic}
Let $\Sigma=\St'(\fm)$ be any costandard module of $G_n$. Then
\begin{equation}\label{eq:acyclic}
 \cE^j_n(\Sigma)\simeq
 \begin{cases}\Sigma&j=0,\\0&j>0.\end{cases}
\end{equation}
In particular, all terms of $\bD^\bullet(\fl)$ are
$\cE^0_n$-acyclic.
\end{lemma}
\begin{proof}
Unlinked segment factors commute
\cite[Theorem~9.7(a)]{Zelevinsky1980},
\cite[\S3.1 and Theorem~3.8]{RegularLM2018}.
Move all factors off the cuspidal line
$\{\nu^{1/2+j}:j\in\ZZ\}$ to the left and remove them using
\eqref{gl:segment-pullout} with $t=0$.
On the remaining line, arrange the factors
$\De_{\one_1}[b_i,a_i]$ by increasing lower endpoint $a_i$.

Continue applying \eqref{gl:segment-pullout} from the left until
the first factor with lower endpoint $1/2$. If no such factor occurs,
the calculation ends at $\cE^0_0(\one_0)=\one_0$, with no higher
Ext groups. Otherwise the remaining product $\Pi$, of rank $N$,
has cuspidal support contained in $\{\nu^u:u\geq1/2\}$.
For $r<N$, the $G_{N-r}$-character in the rank-$r$ piece of
$\omega_{N,N}$ has support
\[
 \nu^{1/2-(N-r)}+\cdots+\nu^{-1/2}.
\]
These characters cannot occur in the cuspidal support of any
Jacquet constituent of $\Pi$. Second adjointness and the K\"unneth argument
of \cite[\S3.3]{CLLTZ2025} therefore give
$\Ext^j_{G_N}(\Xi_r,\Pi)_\sm=0$ for all $j$.
The rank filtration and \eqref{eq:open-piece} give
$\cE^0_N(\Pi)=\Pi$, with no higher Ext groups.
\end{proof}

\section{Termwise theta lifts and induced maps}

We recall the Godement--Jacquet zeta integral \cite[\S\S3--5]{GJ1972}.
To use the same variable as the character-line endpoints, put
\[
 L(s,\Pi)=L_{\mathrm{GJ}}\bigl(es-(e-1)/2,\Pi\bigr).
\]
This is the usual normalization when $D=F$, and $s=1/2$ is the
theta point for every $D$ \cite[Theorem~2.1]{FSX2018}.
Fix a Haar measure $dg$ on $G_n$. For a finite-length representation $\Pi$,
$v\in\Pi$, $\lambda\in\Pi^\vee$, and $\varphi\in\cS(M_{n,n})$, it is
\begin{equation}\label{eq:zeta}
 Z(s,\varphi,v,\lambda)
 =\int_{G_n}\varphi(g)\langle\Pi(g)v,\lambda\rangle
                   \nu_n(g)^{s+(n-1)/2}\,dg.
\end{equation}
It converges in a right half-plane and extends rationally in $q^{-es}$.
For irreducible $\Pi$ the integrals generate the fractional ideal
$L(s,\Pi)\CC[q^{es},q^{-es}]$.

For a standard or costandard module
$\Sigma=\delta_1\times\cdots\times\delta_r$, the generator is
$\prod_iL(s,\delta_i)$.
At $s=1/2$, they define a $G_n$-equivariant map
\begin{equation}\label{eq:zeta-map}
 \cZ_\Pi:\Pi\longrightarrow\cE^0_n(\Pi),\qquad
 \big\langle\cZ_\Pi(v)(\varphi),\lambda\big\rangle
 =\left.\frac{Z(s,\varphi,v,\lambda)}{L(s,\Pi)}\right|_{s=1/2},
\end{equation}
see \cite{FSX2018}. For a segment $[b,a]_\chi$, one has
\[
 L(s,\De_\chi[b,a])=L_F(e(s+b),\chi),\qquad
 L_F(u,\chi)=
 \begin{cases}
 (1-\chi(\varpi)q^{-u})^{-1}&\chi\text{ unramified},\\
 1&\chi\text{ ramified}.
 \end{cases}
\]
This is the discrete-series factor formula
\cite[Theorem~2.2(b),(f),(h) and equation~(8) in its proof]{ABPSDepth2016}
in our variable.
An empty segment contributes $1$. For a ladder $\fl$, multiplicativity
gives
\begin{equation}\label{eq:term-L}
 L_w(s):=L(s,\Upsilon_w(\fl))
 =\prod_{\substack{1\leq r\leq k\\a_r\leq b_{w(r)}}}
   L_F(e(s+b_{w(r)}),\chi)\qquad(w\in W(\fl)).
\end{equation}
\begin{lemma}\label[lemma]{lem:zeta-isomorphism}
For every costandard module $\Sigma$, the map $\cZ_{\Sigma}$ is an isomorphism.
\end{lemma}
\begin{proof}
By \cref{lem:acyclic}, the source and target are isomorphic.
The simple socle of $\Sigma$ occurs with multiplicity one in $\Sigma$, hence any non-zero endomorphism of $\Sigma$ must be an isomorphism.
The normalized zeta integrals are not all zero at $1/2$, hence
$\cZ_{\Sigma}$ is an isomorphism.
\end{proof}

\begin{prop}\label[prop]{prop:maps}
Let $v\lessdot w$ in $W(\fl)$, and let
$p_{v,w}:\Upsilon_v(\fl)\sra \Upsilon_w(\fl)$ be a component of the costandard BGG
differential. Under the zeta isomorphisms,
\begin{equation}\label{eq:map-L-ratio}
 \cE^0_n(p_{v,w})\,\cZ_{\Upsilon_v(\fl)}
 =\left.\frac{L_w(s)}{L_v(s)}\right|_{s=1/2}
   \cZ_{\Upsilon_w(\fl)}\,p_{v,w}.
\end{equation}
The ratio is holomorphic at $1/2$. In particular, the induced map is
nonzero if and only if the pole orders of $L_v$ and $L_w$ at $1/2$ agree.
\end{prop}
\begin{proof}
Functoriality of matrix coefficients gives the identity of unnormalized
zeta integrals. Dividing by $L_v(s)$ and specializing gives
\eqref{eq:map-L-ratio}. Since $\Upsilon_w(\fl)$ is a quotient of $\Upsilon_v(\fl)$,
its zeta ideal is contained in that of $\Upsilon_v(\fl)$, so $L_w/L_v$ is a
Laurent polynomial in $q^{-es}$. The last assertion follows from
\eqref{eq:term-L} and \cref{lem:zeta-isomorphism}.
\end{proof}

\begin{corollary}\label[corollary]{cor:poles}
For $w\in W(\fl)$, the function $L_w(s)$ has at most a simple pole at
$s=1/2$. It has a pole if and only if $\fl_w$ contains a nonempty segment
with upper endpoint $\nu^{-1/2}$.
\end{corollary}
\section{Equal rank for balanced multisegments}\label{sec:regular-equal-rank}

We compute the equal-rank higher theta lifts for the balanced
multisegments of Lapid--M\'inguez \cite[\S6]{RegularLM2018}.
The ladder formula follows in the next section.

\subsection{Balanced multisegments and their resolutions}
\label{gl:regular:resolution-subsection}

Fix a character $\chi$ of $\Ff^\times$. A multisegment
\[
 \fm=[b_1,a_1]_\chi+\cdots+[b_k,a_k]_\chi
\]
on its cuspidal line is \emph{regular} if the lower endpoints
$a_1,\ldots,a_k$ are pairwise distinct and the upper endpoints
$b_1,\ldots,b_k$ are pairwise distinct. The segments are nonempty,
and all endpoint differences are integral.

We define balancedness by the following pattern avoidance condition.
Choose a sub-multisegment consisting of $r\geq4$ segments and write it
as $\sum_{i=1}^r[d_i,c_i]_\chi$, ordered so that
$d_1>\cdots>d_r$. It is of \emph{type $4231$} if
\[
 \begin{aligned}
  c_{i+1}<c_i\leq d_{i+1}+1&\qquad(3\leq i<r),\\
  c_3<c_1\leq d_3+1,&\\
  c_r<c_2<c_{r-1}.&
 \end{aligned}
\]
It is of \emph{type $3412$} if
\[
 \begin{aligned}
  c_{i+1}<c_i\leq d_{i+1}+1&\qquad(4\leq i<r),\\
  c_4<c_2\leq d_4+1,&\\
  c_3<c_r<c_1<c_\ell,&
 \end{aligned}
 \qquad
 \ell=\begin{cases}2&r=4,\\r-1&r>4.\end{cases}
\]
The regular multisegment $\fm$ is \emph{balanced} if it contains
no sub-multisegment of either type. This is the characterization of
Lapid--M\'inguez in \cite[Definition~6.10 and Proposition~6.13]{RegularLM2018}.

These conditions can be expressed geometrically as follows. Order the lower endpoints
increasingly as $a_1<\cdots<a_k$ and the upper endpoints
decreasingly as $b_1>\cdots>b_k$. For $\tau\in S_k$, define
\begin{equation}\label{gl:regular:data}
	\fm_\tau=\sum_{i=1}^k[b_{\tau(i)},a_i]_\chi.
\end{equation}
As before, $[a-1,a]_\chi$ contributes
$\one_0$, while a pair with $a>b+1$ makes the whole induced term zero. There is a unique
$\sigma\in S_k$ such that the segment with lower endpoint $a_i$
has upper endpoint $b_{\sigma(i)}$. Thus $\fm=\fm_\sigma$. Unlike $\fm_\sigma$,
an alternative pairing may contain empty or inadmissible segments.
All permutation orders below are Bruhat orders.
There is a unique minimal admissible permutation $\sigma_0$ such that
\begin{equation}\label{gl:regular:admissible}
 a_i\leq b_{\tau(i)}+1\text{ for every }i
 \quad\Longleftrightarrow\quad \sigma_0\leq\tau.
\end{equation}
See \cite[\S6.1, equation~(15)]{RegularLM2018}.

Let $C_w$ be the Schubert cell indexed by $w$ in the full flag variety
of $\GL_k(\CC)$, and let $X_w=\overline{C_w}$.
It is known that the pattern avoidance condition above
is equivalent to $X_\sigma$ being smooth along $C_{\sigma_0}$,
and also to its being rationally smooth there
\cite[Corollary~6.7 and Theorem~1.2(4), (6)]{RegularLM2018}.
The underlying equivalence of smoothness and rational smoothness in
type $A$ is due to Deodhar \cite{Deodhar1985}.
Equivalently, the Kazhdan--Lusztig polynomials satisfy
\begin{equation}\label{gl:regular:smooth}
 P_{\tau,\sigma}(q)=1\qquad(\sigma_0\leq\tau\leq\sigma).
\end{equation}
Throughout this section we assume that $\fm_\sigma$ is balanced.

Put $\pi=\Lrep(\fm_\sigma)$ and define as for ladder representations
\begin{align}
 \Sigma_\tau&=\De_\chi[b_1,a_{\tau^{-1}(1)}]\times\cdots\times
              \De_\chi[b_k,a_{\tau^{-1}(k)}],\notag\\
 \Upsilon_\tau&=\De_\chi[b_{\tau(1)},a_1]\times\cdots\times
              \De_\chi[b_{\tau(k)},a_k].
 \label{gl:regular:terms}
\end{align}
\begin{theorem}\label[theorem]{gl:regular:resolution}
Let $\fm_\sigma$ be balanced of length $N$, and put
$R=\ell(\sigma)-\ell(\sigma_0)$. There are exact sequences
\begin{align}
 0\longrightarrow C_R\longrightarrow\cdots\longrightarrow C_0
   \longrightarrow\pi\longrightarrow0,
   \label{gl:regular:standard-resolution}\\
 0\longrightarrow\pi\longrightarrow D^0\longrightarrow\cdots
   \longrightarrow D^R\longrightarrow0,
   \label{gl:regular:costandard-coresolution}
\end{align}
where
\begin{equation}\label{gl:regular:resolution-terms}
 C_r=\bigoplus_{\substack{\sigma_0\leq\tau\leq\sigma\\
                         \ell(\sigma)-\ell(\tau)=r}}\Sigma_\tau,
 \qquad
 D^r=\bigoplus_{\substack{\sigma_0\leq\tau\leq\sigma\\
                         \ell(\sigma)-\ell(\tau)=r}}\Upsilon_\tau.
\end{equation}
The component $\Sigma_\tau\to \Sigma_\upsilon$ on a cover
$\tau\lessdot\upsilon$ is the image of a coherently signed Verma
embedding and is a nonzero injection. The corresponding component
$\Upsilon_\upsilon\to \Upsilon_\tau$ is a nonzero surjection. All multiplicities are
one, and these actual differentials square to zero.
\end{theorem}

\begin{proof}
We first construct the standard resolution and then obtain the
costandard coresolution by duality.

For the upper triangular Borel of $\mathfrak{gl}_k$, the half sum of
positive roots is
\[
 \rho_k=\left(\frac{k-1}{2},\frac{k-3}{2},\ldots,\frac{1-k}{2}\right).
\]
Set
\[
 \lambda=(b_1+1,\ldots,b_k+1)-\rho_k,\qquad
 \mu_\tau=\tau(a_1,\ldots,a_k)-\rho_k.
\]
Use the block $\mathcal O_{\mu_e}$ and the standard, simple and
costandard objects of \cref{prelim:category-o}, with highest weights
$\mu_\tau$.
The weight $\mu_e$ is antidominant because its shifted simple-root
pairings are $a_i-a_{i+1}<0$. With this indexing, a Bruhat cover
$\tau\lessdot\upsilon$ gives a Verma embedding
$\Sigma_\tau^{\mathcal O}\hra \Sigma_\upsilon^{\mathcal O}$
\cite[\S1, equation~(1.2)]{RegularDhillon2019}.

Put $V=(\CC^k)^{\otimes N}$. Since $\lambda+\rho_k$ is strictly
dominant, \cref{prelim:arakawa-suzuki} gives an exact functor
\[
 \mathcal F:\mathcal O_{\mu_e}\longrightarrow\rep_N
\]
whose image has cuspidal support on the $\chi$-line.
The weight of $V$ required for the image of $\Sigma_\tau^{\mathcal O}$ is
\[
 \lambda-\mu_\tau
   =(b_i+1-a_{\tau^{-1}(i)})_{i=1}^k.
\]
These coordinates sum to $N$. They form a weight of $V$ exactly
when they are nonnegative integers, which is equivalent to
$\tau\geq\sigma_0$. The $i$th coordinate is the length of the factor
$\De_\chi[b_i,a_{\tau^{-1}(i)}]$ in $\Sigma_\tau$.
If it is zero, that factor is the empty segment, contributing
$\one_0$. Hence \cite[Theorem~4.2.3]{BC2015} gives
\[
 \mathcal F(\Sigma_\tau^{\mathcal O})=
 \begin{cases}
  \Sigma_\tau&\tau\geq\sigma_0,\\
  0&\tau\not\geq\sigma_0,
 \end{cases}
 \qquad
 \mathcal F(\pi_\sigma^{\mathcal O})=\pi.
\]
The simple-image condition in that theorem is automatic because
$\lambda+\rho_k$ is strictly dominant.

We now prove exactness of the proposed resolution. Work in the Serre subcategory generated by
the $\pi_\tau^{\mathcal O}$ with $\tau\leq\sigma$, whose maximum
length is $\ell(\sigma)$. Dhillon's bounded spectral sequence
\cite[Corollaries~2.9 and 4.3]{RegularDhillon2019} gives
\begin{equation}\label{gl:regular:dhillon}
 E_1^{p,q}=
 \bigoplus_{\substack{\tau\leq\sigma\\
                     \ell(\tau)=\ell(\sigma)+p}}
 \Sigma_\tau^{\mathcal O}\otimes
 \Ext^{-p-q}_{\mathcal O}
       (\pi_\sigma^{\mathcal O},\Upsilon_\tau^{\mathcal O})^*
 \ \Longrightarrow\ \mathrm H^{p+q}(\pi_\sigma^{\mathcal O}).
\end{equation}
Here $\pi_\sigma^{\mathcal O}$ is regarded as a complex concentrated in degree $0$.
The description following Corollary~2.9 identifies the $d_1$ maps in
that row with the coherently signed Verma embeddings.
For $\sigma_0\leq\tau\leq\sigma$, \eqref{gl:regular:smooth}, the
symmetry $P_{\tau^{-1},\sigma^{-1}}=P_{\tau,\sigma}$ and the
intersection-cohomology stalk formula of Kazhdan--Lusztig
\cite{KL1980} (see also \cite[proof of Theorem~3.11.4(iv), p.~512]{BGS1996}) give
\[
 j_{\tau^{-1}}^*\operatorname{IC}(X_{\sigma^{-1}})
 =\CC_{C_{\tau^{-1}}}[\ell(\sigma)].
\]
The dictionary of \cref{prelim:category-o}, adjunction and the
contractibility of $C_{\tau^{-1}}$ give
\[
 R\operatorname{Hom}_{\mathcal O}
 (\pi_\sigma^{\mathcal O},\Upsilon_\tau^{\mathcal O})
 \simeq R\Gamma(C_{\tau^{-1}},\CC)[\ell(\tau)-\ell(\sigma)]
 \simeq\CC[\ell(\tau)-\ell(\sigma)].
\]
This is the geometric Ext computation of
\cite[Remark~4.6]{RegularDhillon2019}, so
\[
 \Ext^j_{\mathcal O}(\pi_\sigma^{\mathcal O},\Upsilon_\tau^{\mathcal O})
 \simeq
 \begin{cases}
  \CC&j=\ell(\sigma)-\ell(\tau),\\
  0&\text{otherwise}.
 \end{cases}
\]
Apply $\mathcal F$ to the entire spectral sequence. Exactness
preserves its pages and abutment. Terms with $\tau\not\geq\sigma_0$
vanish, and the displayed Ext computation puts every surviving term
in the row $q=0$. More precisely,
\[
 \mathcal F(E_1^{-r,0})=C_r,\qquad
 \mathcal F(E_1^{p,q})=0\quad(q\ne0).
\]
Thus this row is the complex
$\cdots\to C_1\to C_0$, with $C_r$ in cohomological degree $-r$.
There are no higher differentials, so its cohomology is the abutment,
namely $\pi$ in degree zero and zero elsewhere.
The image of $\Sigma_\sigma^{\mathcal O}\to \pi_\sigma^{\mathcal O}$
supplies the augmentation. This proves
\eqref{gl:regular:standard-resolution}.
Each cover map is a transported Verma embedding. Exactness makes it
injective, and its nonzero domain makes it nonzero.

For the coresolution, replace the endpoint data by
$a_i^\vee=-b_i$, $b_i^\vee=-a_i$ and $\chi^\vee=\chi^{-1}$.
Writing the endpoint data explicitly, one has
$\fm_{\tau^{-1}}(a^\vee,b^\vee)=\fm_\tau(a,b)^\vee$.
Inversion preserves Bruhat order, length and Kazhdan--Lusztig
polynomials, so the dual data are balanced, with minimal admissible
permutation $\sigma_0^{-1}$ and indexing permutation $\sigma^{-1}$.
Apply the standard resolution just proved to these dual data and
take smooth contragredients. Since
\[
 \bigl(\Sigma_{\tau^{-1}}(a^\vee,b^\vee)\bigr)^\vee
 \simeq\Upsilon_\tau(a,b),
\]
exactness of smooth duality gives
\eqref{gl:regular:costandard-coresolution} with the stated terms.
Its cover components are the duals of the corresponding nonzero
injections, hence are nonzero surjections, and their differentials
square to zero.
\end{proof}

The construction also works for balanced multisegments on any fixed supercuspidal line
$\rho$ of $\GL_h(D)$, using its step $\nu_\rho$.
\begin{prop}\label[prop]{gl:regular:theta-complex}
Assume that the cuspidal line is that of
$\chi\circ\operatorname{Nrd}$ on $D^\times$, where $\chi$ is a
smooth character of $\Ff^\times$, and let $\pi=\Lrep(\fm_\sigma)$
be as above. Then
\begin{equation}\label{gl:regular:theta-complex-equation}
 \cE_N^r(\pi)\simeq
 \mathrm H^r\bigl(\cE_N^0(D^\bullet)\bigr).
\end{equation}
For every signed component $p_{\upsilon,\tau}:\Upsilon_\upsilon\to \Upsilon_\tau$
of the coresolution differential, where $\tau\lessdot\upsilon$,
the zeta isomorphisms satisfy
\begin{equation}\label{eq:regular-map-L-ratio}
 \cE_N^0(p_{\upsilon,\tau})\,\cZ_{\Upsilon_\upsilon}
 =\left.\frac{L(s,\Upsilon_\tau)}{L(s,\Upsilon_\upsilon)}\right|_{s=1/2}
   \cZ_{\Upsilon_\tau}\,p_{\upsilon,\tau}.
\end{equation}
The ratio is holomorphic at $1/2$.
\end{prop}

\begin{proof}
Every term of the coresolution is acyclic by \cref{lem:acyclic}, and
\cref{lem:acyclic-coresolution} gives
\eqref{gl:regular:theta-complex-equation}. The proof of
\cref{prop:maps} applies to every costandard quotient map:
functoriality of matrix coefficients gives
\eqref{eq:regular-map-L-ratio}, and containment of zeta ideals
gives holomorphy of the ratio.
\end{proof}

\subsection{Deletion at an admissibility corner}

Let $D^\bullet(\fm_\sigma)$ be the costandard coresolution of
\cref{gl:regular:resolution}. Assume there exist indices $i,j$ with $a_i=b_j+1$.
We consider the summands indexed by
\[
 J=\{\tau\in[\sigma_0,\sigma]:\tau(i)=j\}.
\]
For these indices, the $i$th factor of $\Upsilon_\tau$ is
$\De_\chi[b_j,a_i]=\one_0$. In each degree put
\[
 B^r=\bigoplus_{\substack{\tau\in J\\
                         \ell(\sigma)-\ell(\tau)=r}}\Upsilon_\tau
 \ \subseteq\ D^r(\fm_\sigma).
\]
We will show that the differential of $D^\bullet(\fm_\sigma)$
preserves these summands and identify the resulting subcomplex.

Delete $a_i$ and $b_j$ from the endpoint lists, retaining their
orders, and call the shorter lists $a',b'$. If $\tau(i)=j$, let
$\operatorname{del}(\tau)\in S_{k-1}$ be obtained by deleting the
entry $j$ at position $i$ and renumbering the remaining values
increasingly. Write $\operatorname{ins}$ for the inverse operation
and $\sigma'_0$ for the minimal admissible permutation for $a',b'$.

\begin{lemma}\label[lemma]{gl:regular:corner}
The set $J$ is nonempty and has a unique maximum $\sigma_*$.
Put $\sigma'=\operatorname{del}(\sigma_*)$. Deletion gives an
order isomorphism
\[
 \operatorname{del}:J\xrightarrow{\sim}[\sigma'_0,\sigma'].
\]
The multisegment
\begin{equation}\label{eq:regular-star}
 \fm^*=\sum_{r\ne i}[b_{\sigma_*(r)},a_r]_\chi
\end{equation}
has only nonempty segments, has the same total length as $\fm_\sigma$,
and is balanced. Put
\begin{equation}\label{eq:regular-degree}
 d=\ell(\sigma)-\ell(\sigma_*)
   =\ell(\sigma)-(2k-i-j)-\ell(\sigma').
\end{equation}
Then $d\geq1$. The differentials of $D^\bullet(\fm_\sigma)$
restrict to maps $B^r\to B^{r+1}$, and the resulting complex satisfies
\begin{equation}\label{eq:regular-corner-complex}
 B^\bullet\simeq D^\bullet(\fm^*)[-d].
\end{equation}
In particular, $H^d(B^\bullet)=\Lrep(\fm^*)$ and all its other
cohomology groups vanish.
\end{lemma}

\begin{proof}
First, $J$ is nonempty. Put $p=\sigma^{-1}(j)$ and $q=\sigma(i)$.
Since the segments of $\fm_\sigma$ are nonempty, $p<i$ and $q<j$.
Swap the values of $\sigma$ at $p$ and $i$, obtaining
$\tau_1<\sigma$ with $\tau_1(i)=j$. Its other segments are
nonempty, so $\tau_1\in J$.

For any admissible $\tau$ with $\tau(i)=j$, the endpoint orders give
$r>i\Rightarrow\tau(r)<j$ and $\tau(r)>j\Rightarrow r<i$.
Thus the entry $(i,j)$ contributes exactly $(k-i)+(k-j)$ inversions.
Deletion preserves Bruhat order by
\cite[Lemma~17]{BilleyWarrington2003}, recalled in
\cite[Lemma~5.5]{RegularLM2018}, and the preceding inversion count
shows that it subtracts $2k-i-j$ from the length. It also identifies admissibility
with that for $a',b'$.

Put $\tau_0=\operatorname{ins}(\sigma'_0)$. Then
$\sigma_0\leq\tau_0\leq\tau_1\leq\sigma$.
The subgroup $\{w\in S_k:w(j)=j\}\simeq S_{k-1}$ acts by
permuting the values other than $j$. Its left coset through $\tau_0$
consists of the permutations $\tau$ with $\tau(i)=j$.
On this coset, the map used in \cite[\S2.3]{RegularBilleyBraden2003}
is $\operatorname{del}$: it records the relative order of the entries
remaining after $j$ is removed.
Balancedness gives $P_{\tau_0,\sigma}(1)=1$, and
\cite[Theorem~4]{RegularBilleyBraden2003} gives
\[
 1\geq\sum_{y\in\mathcal M}
       P_{y,\sigma}(1)
       P'_{\sigma'_0,\operatorname{del}(y)}(1),
\]
where $\mathcal M$ is the set of maximal elements of
$\{y\leq\sigma:y(i)=j\}$, and $P'$ denotes the Kazhdan--Lusztig
polynomials for $S_{k-1}$. Exactly the maxima above $\tau_0$
contribute, each by a positive integer. At least one exists, so there
is exactly one, denoted $\sigma_*$, and
$P'_{\sigma'_0,\sigma'}=1$. Every element of $J$ lies between
$\tau_0$ and this maximum. Deletion therefore identifies $J$ with
$[\sigma'_0,\sigma']$.

By \eqref{gl:regular:admissible} with upper endpoints $b'_s-1$,
the inequalities $a'_r\leq b'_{u(r)}$ for all $r$ remain valid
when $u$ is replaced by a larger permutation in Bruhat order.
They hold for $u=\operatorname{del}(\tau_1)$, hence also for
$\sigma'\geq\operatorname{del}(\tau_1)$. Thus all segments of
$\fm^*$ are nonempty. The identity
$P'_{\sigma'_0,\sigma'}=1$ makes it balanced by
\cite[\S5.1 and Corollary~6.7]{RegularLM2018}.
Deleting an empty segment leaves the total length unchanged.
Moreover, $\sigma_*<\sigma$, so $d\geq1$, and the length calculation
above gives \eqref{eq:regular-degree}.

Finally, a component of the differential out of $\Upsilon_\tau$ has target
$\Upsilon_\eta$ with $\eta\lessdot\tau$. If $\tau(i)=j$ and
$\eta(i)\ne j$, the transposition creates either an entry $(i,>j)$
or an entry $(>i,j)$. Both violate admissibility. Hence the differential
sends $B^r$ into $B^{r+1}$.
Omitting the empty factor identifies $\Upsilon_\tau$ with the smaller
costandard indexed by $\operatorname{del}(\tau)$, and
\[
 \ell(\sigma)-\ell(\tau)
 =d+\ell(\sigma')-\ell(\operatorname{del}(\tau)).
\]
To identify the maps, dualize and normalize the standard embeddings
into $\Upsilon_{\sigma_*}^\vee$. The one-dimensional Hom spaces of
\cref{gl:line-transport} identify them with the smaller standard
embeddings. Deletion preserves covers and length-two diamonds,
so choose the same coherent BGG signs. Dualizing gives
\eqref{eq:regular-corner-complex}, and
\cref{gl:regular:resolution} gives its asserted cohomology.
\end{proof}

\subsection{The higher theta representations}

\begin{theorem}\label[theorem]{gl:regular:theta}
Let $\fm_\sigma$ be a balanced multisegment on a character line,
with all its segments nonempty, and put $\pi=\Lrep(\fm_\sigma)$.
If its lower endpoint characters do not include $\nu^{1/2}$, or its upper endpoint
characters do not include $\nu^{-1/2}$, then
\[
 \Theta_N^0(\pi)\simeq\pi^\vee,
 \qquad \Theta_N^r(\pi)=0\quad(r>0).
\]
Otherwise let $i,j$ be the unique indices such that
$\chi\nu^{a_i}=\nu^{1/2}$ and
$\chi\nu^{b_j}=\nu^{-1/2}$. Define $\fm^*,d$ by
\eqref{eq:regular-star}--\eqref{eq:regular-degree} and put
$\pi^*=\Lrep(\fm^*)$. Then
\begin{equation}\label{eq:regular-theta}
 \Theta_N^r(\pi)\simeq
 \begin{cases}
  \pi^\vee&r=0,\ d>1,\\
  (\pi^*)^\vee&r=d,\\
  (\pi^*)^\vee&r=d-1,\ d>1,\\
  \Sigma(\pi)&r=0,\ d=1,\\
  0&\text{otherwise},
 \end{cases}
\end{equation}
where $\Sigma(\pi)$ is the nonsplit extension
\begin{equation}\label{eq:regular-extension}
 0\longrightarrow(\pi^*)^\vee\longrightarrow\Sigma(\pi)
   \longrightarrow\pi^\vee\longrightarrow0.
\end{equation}
\end{theorem}

\begin{proof}
The costandard coresolution of \cref{gl:regular:resolution} and
the termwise calculation of \cref{gl:regular:theta-complex} give
\[
 \cE_N^r(\pi)
  \simeq H^r\bigl(\cE_N^0(D^\bullet)\bigr).
\]
By multiplicativity of the segment $L$-factors, the term indexed by $\tau$ has a
simple pole exactly when the upper endpoint character $\nu^{-1/2}$ occurs in a
nonempty segment. When there is no critical upper endpoint, every pole order
is zero. When there is a critical upper endpoint but no critical lower endpoint, every
term has the same simple pole. For each term set
\[
 c_\tau=\left.(1-q^{-e(s-1/2)})^{\epsilon_\tau}
                 L(s,\Upsilon_\tau)\right|_{s=1/2}\ne0,
\]
where $\epsilon_\tau$ is its pole order, and rescale its zeta
identification by $c_\tau$.
Proposition~\ref{gl:regular:theta-complex} then leaves the
original differential when the pole orders agree and gives zero
otherwise. In the two constant-pole cases this identifies the
theta complex with $D^\bullet$ and proves the first assertion.

In the critical case, a term has no pole precisely when $\tau(i)=j$.
After this rescaling, differentials within either pole class
are the original differentials, and all differentials between the
classes vanish. Therefore there is an isomorphism of actual complexes
\[
 \cE_N^0(D^\bullet)\simeq (D^\bullet/B^\bullet)\oplus B^\bullet.
\]
By \cref{gl:regular:corner}, $B$ has just $\pi^*$ in degree $d$.
The long exact sequence of $0\to B\to D\to D/B\to0$ now gives,
for $d>1$, $H^0(D/B)=\pi$, $H^{d-1}(D/B)=\pi^*$, and no other
cohomology. For $d=1$ it gives
\[
 0\longrightarrow\pi\longrightarrow H^0(D/B)
    \longrightarrow\pi^*\longrightarrow0,
 \qquad H^r(D/B)=0\quad(r>0).
\]
The two irreducibles are distinct, because their Langlands
multisegments have respectively $k$ and $k-1$ nonempty segments.
Howe duality \cite[Theorem~1 and Corollary~6.3]{Minguez2008} says
that $\Theta_N^0(\pi)$ has the unique irreducible quotient
$\pi^\vee$. Hence the last extension cannot split. Duality
gives \eqref{eq:regular-theta}--\eqref{eq:regular-extension}.
\end{proof}

\section{Higher theta lifts of ladders}

We make the results of the previous section more explicit for ladders. 

\subsection{The BGG resolution and coresolution}\label{gl:line-transport}

We take $\rho=\chi\circ\operatorname{Nrd}$ for a character $\chi$.
We may absorb a common shift of the endpoints into $\chi$.

Let $\fl$ be a ladder with $k$ segments. We use the ordinary Bruhat order
on $S_k$ and write $v\lessdot w$ if $v<w$ and $\ell(w)=\ell(v)+1$.
For $v\leq w$ in $W(\fl)$, one has injections
\[
 \iota_{w,v}:\Sigma_w(\fl)\hra \Sigma_v(\fl),
\]
obtained from the maps between Verma modules by the exact functor to
$p$-adic representations recalled in \cref{prelim:arakawa-suzuki}.
The relevant Hom spaces are one-dimensional in the split model
\cite[Proposition~4.3.6]{EmertonHelm2014}, and the equivalence
preserves this property over $D$. We normalize the Verma embeddings
before transport so that their compositions agree.
The set $W(\fl)$ is
a lower Bruhat interval \cite[Lemma~15]{LM2014}.

Choose signs $\varepsilon(w,v)\in\{1,-1\}$ on the Bruhat covers so that, for
each length-two interval $[v,w]$ with intermediate elements $s,t$,
\[
 \varepsilon(w,s)\varepsilon(s,v)
 =-\varepsilon(w,t)\varepsilon(t,v).
\]
We define the BGG complex $\bC^\bullet(\fl)$ by
\[
 \bC^{-j}(\fl)=\bigoplus_{\substack{w\in S_k\\\ell(w)=j}}\Sigma_w(\fl),
\]
with differential $d_C^{-j}:\bC^{-j}\longrightarrow\bC^{-j+1}$ whose
$(w,v)$ component is $\varepsilon(w,v)\iota_{w,v}$ for $v\lessdot w$
and zero otherwise. Similarly, define
\[
 \bD^j(\fl)=\bigl(\bC^{-j}(\fl^\vee)\bigr)^\vee
 =\bigoplus_{\substack{w\in W(\fl)\\\ell(w)=j}}\Upsilon_w(\fl),
\]
where the last equality uses the corresponding relabeling of permutations.
Its differential components are signed surjections
$p_{v,w}:\Upsilon_v(\fl)\sra \Upsilon_w(\fl)$ for $v\lessdot w$.

\begin{theorem}[BGG resolution]\label{thm:bgg}
For any ladder $\fl$, one has
\[
 \mathrm{H}^j(\bC^\bullet(\fl))\simeq \mathrm{H}^j(\bD^\bullet(\fl))\simeq
 \begin{cases}\Lrep(\fl)&j=0,\\0&j\ne0.\end{cases}
\]
Thus $\bC^\bullet(\fl)$ is a resolution and $\bD^\bullet(\fl)$ is a
coresolution of $\Lrep(\fl)$.
\end{theorem}
\begin{proof}
Apply the exact equivalence above to the split BGG resolution
obtained from \cite[Theorem~5.1.1]{Suzuki1998}, in the ladder
conventions of \cite[\S4]{BC2015}. This gives its terms, differentials and
augmentation over $D$. Apply the same construction to the dual
endpoint data and then take the smooth contragredient to obtain the
coresolution.
\end{proof}

\subsection{The higher theta lifts}\label{sec:unequal}
Let $m=n+t$ with $t\geq0$, and put $s_t=(t+1)/2$.
Let $\fl=\sum_{h=1}^k[b_h,a_h]_\chi$ be a ladder of length $n$.
If it has a lower endpoint $\nu^{s_t}$ and an upper endpoint
$\nu^{-s_t}$, let $i<j$ be the unique indices such that
\begin{equation}\label{eq:unequal-critical}
 \chi\nu^{a_i}=\nu^{s_t},\qquad
 \chi\nu^{b_j}=\nu^{-s_t}.
\end{equation}
In this case define
\begin{equation}\label{eq:unequal-star}
 \fl_t^*=\sum_{h<i}[b_h,a_h]_\chi
       +\sum_{h=i}^{j-1}[b_h,a_{h+1}]_\chi
       +\sum_{h>j}[b_h,a_h]_\chi,
 \qquad d_t(\fl)=j-i-t.
\end{equation}
Thus $\fl_t^*$ is obtained by deleting $a_i$ and $b_j$ and pairing
the remaining endpoints in order. It is a ladder of length $n+t=m$.
If either critical endpoint is absent, put $d_t(\fl)=0$.
For $t>0$, the integer $d_t(\fl)$ need not be positive even when
both endpoints occur. In equal rank we write $d(\fl)=d_0(\fl)$
and $\fl^*=\fl_0^*$ when the critical endpoints occur, with
$\fl^*=\fl$ otherwise.

\Needspace{13\baselineskip}
\begin{theorem}\label[theorem]{thm:main}\label[theorem]{thm:unequal}
Let $\fl$ be a ladder of length $n$, let $m=n+t$ with $t\geq0$,
and put $d=d_t(\fl)$. For $r>0$, one has
\begin{equation}\label{eq:unequal-positive}
 \Theta_m^r(\Lrep(\fl))\simeq
 \begin{cases}
  \Lrep(\fl_t^*)^\vee&r=d,\ d\geq1,\\
  \Lrep(\fl_t^*)^\vee&r=d-1,\ d>1,\\
  0&\text{otherwise}.
 \end{cases}
\end{equation}
In degree zero, the following equality holds in the Grothendieck group:
\begin{equation}\label{eq:unequal-zero}
 [\Theta_m^0(\Lrep(\fl))]
 =[\Lrep(\fl)^\vee\times\one_t]
 +\begin{cases}
   [\Lrep(\fl_t^*)^\vee]&d=1,\\
   0&d\ne1.
  \end{cases}
\end{equation}
For $t=0$, one has
$\Theta_n^0(\Lrep(\fl))\simeq\Lrep(\fl)^\vee$ if $d\ne1$.
If $t=0$ and $d=1$, the lift is the nonsplit extension
$\Sigma(\fl)$ given by
\begin{equation}\label{eq:Sigma}
 0\longrightarrow\Lrep(\fl_0^*)^\vee
 \longrightarrow\Sigma(\fl)
 \longrightarrow\Lrep(\fl)^\vee\longrightarrow0.
\end{equation}
\end{theorem}

We first recall the Jacquet-component functors of
\cite[\S3]{Minguez2009}, see \cite[\S\S5--6]{Minguez2009}
for the associated socle and combinatorial results. For a finite-length
representation $\Pi$ of $G_N$, let $\partial_x\Pi$ be the generalized
$\nu^x$-isotypic component of $r_{P_{1,N-1}}\Pi$, viewed as a
$G_{N-1}$-representation. In particular, $\partial_x$ is exact.
For a ladder $\fl=\sum_h[b_h,a_h]$ on the line of real powers of
$\nu$, the Jacquet module is semisimple. If $x=b_h$
and either $h=k$ or $b_h>b_{h+1}+1$, a standard computation gives
\begin{equation}\label{eq:ladder-derivative}
 \partial_x\Lrep(\fl)\simeq
 \Lrep(\fl-[b_h,a_h]+[b_h-1,a_h]),
\end{equation}
where empty segments are omitted. All other derivatives vanish. See
\cite[Theorem~2.1 and Corollary~2.2]{KL2012}.
When \eqref{eq:unequal-critical} holds, absorb $\chi$ into the
endpoints, so that $a_i=s_t$ and $b_j=-s_t$. Then
\eqref{eq:ladder-derivative} implies
\begin{equation}\label{eq:exceptional-derivative-zero}
 \partial_{1-s_t}\Lrep(\fl)=0.
\end{equation}
Indeed, an upper endpoint $1-s_t$ must be $b_{j-1}$ and cannot be shortened because
$b_j=-s_t$.

\begin{lemma}\label[lemma]{lem:theta-derivative}
Let $\pi$ be a ladder representation whose cuspidal support consists of
characters $\nu^x$ with $x\in\mathbb R$, and let $m>n\geq1$.
For every $r\geq0$ and $x\ne1-s_t$, one has
\begin{equation}\label{eq:ordinary-recursion}
 \partial_x\cE^r_{n,m}(\pi)
 \simeq\cE^r_{n-1,m-1}(\partial_x\pi).
\end{equation}
If $\partial_{1-s_t}\pi=0$, then
\begin{equation}\label{eq:exceptional-recursion}
 \partial_{1-s_t}\cE^r_{n,m}(\pi)
 \simeq\nu^{1/2}\cE^r_{n,m-1}(\nu^{-1/2}\pi).
\end{equation}
\end{lemma}
\begin{proof}
Put $Q=P_{1,m-1}$. We first have the derived Jacquet identity
\begin{equation}\label{eq:derived-jacquet}
 r_Q\cE^r_{n,m}(\pi)
 \simeq\Ext^r_{G_n}(r_{\overline Q}\omega_{n,m},\pi)_\sm.
\end{equation}
To justify this identity, choose a projective resolution
$P_\bullet\to\omega_{n,m}$ in $\rep(G_n\times G_m)$ and put
$T_\bullet=(P_\bullet\otimes\pi^\vee)_{G_n}$.
Restriction to $G_n$ preserves projectives, as one sees using
compact induction from product compact open subgroups. Admissibility
of $\pi$ identifies $\Hom_{G_n}(P_\bullet,\pi)_\sm$ with
$T_\bullet^\vee$. Exactness of smooth contragredients
\cite[Proposition~2.1.11]{Casselman1995} therefore gives
$H_r(T_\bullet)^\vee\simeq\cE^r_{n,m}(\pi)$.
The natural injection
$H_r(T_\bullet)\hra\cE^r_{n,m}(\pi)^\vee$ and
\cref{gl:rank-homology} show that $H_r(T_\bullet)$ has finite length.
Casselman's Jacquet duality \cite[Corollary~4.2.5]{Casselman1995},
together with exactness of $r_{\overline Q}$ and its commutation
with $G_n$-coinvariants, now yields
\[
 \begin{aligned}
 r_Q\cE^r_{n,m}(\pi)
 &\simeq\bigl(r_{\overline Q}H_r(T_\bullet)\bigr)^\vee\\
 &\simeq H_r\bigl((r_{\overline Q}P_\bullet
                   \otimes\pi^\vee)_{G_n}\bigr)^\vee\\
 &\simeq\Ext^r_{G_n}(r_{\overline Q}\omega_{n,m},\pi)_\sm.
 \end{aligned}
\]
For the last isomorphism, $r_{\overline Q}P_\bullet$ is a projective
resolution: the Jacquet functor is left adjoint to exact parabolic
induction, and restriction from $G_n\times M$ to $G_n$ preserves
projectives, where $M$ is the Levi factor of $Q$.

The two-step Kudla filtration of $r_{\overline Q}\omega_{n,m}$
\cite[Proposition~3.3]{Minguez2008}, followed by second adjointness and
the regular representation calculation, gives a long exact sequence
\[
 \cdots\longrightarrow \Gamma^r\longrightarrow r_Q\cE^r_{n,m}(\pi)
 \longrightarrow \Lambda^r\longrightarrow \Gamma^{r+1}\longrightarrow\cdots,
\]
where, by a standard computation with the normalizing characters,
\begin{align*}
 \Gamma^r&=\nu^{1-s_t}\boxtimes
       \nu^{1/2}\cE^r_{n,m-1}(\nu^{-1/2}\pi),\\
 \Lambda^r&=\bigoplus_x\nu^x\boxtimes
       \cE^r_{n-1,m-1}(\partial_x\pi).
\end{align*}
Here we use the semisimplicity of the Jacquet module of $\pi$. Taking
the generalized $\nu^x$-component proves
\eqref{eq:ordinary-recursion}. At $x=1-s_t$, the $\Lambda^r$-terms vanish
under the additional hypothesis, giving
\eqref{eq:exceptional-recursion}.
\end{proof}

\begin{lemma}[Support]\label[lemma]{gl:support}
If $\pi$ is a ladder whose cuspidal support consists of real powers of
$\nu$, the same is true of every constituent of $\cE^i_{n,m}(\pi)$.
\end{lemma}
\begin{proof}
Use the full matrix-rank filtration of $\omega_{n,m}$. The normalized
Jacquet modules of $\pi$ are semisimple direct sums by
\cite[Corollary~2.2]{KL2012} and the transfer in
\cref{gl:line-transport}, with cuspidal factors drawn from those
of $\pi$. The rank-stratum calculation
\cite[\S5.1, equation~(10), Proposition~5.9]{APS2017} expresses each stratum's
smooth Ext as a sum of target parabolic inductions from a real reduced-norm
twist of a source Jacquet factor on $\GL_k(D)$ and a real reduced-norm
character on $G_{m-k}$, tensored with a finite-dimensional vector
space. The latter is a residual-source Ext group, with trivial target
action. Its finiteness is \cite[Proposition~3.6]{APS2017}.
Thus every stratum Ext has the asserted support. The smooth long exact
sequences of \cite[Lemma~5.14]{APS2017} and closure under subquotients
and extensions propagate it through the finite filtration. This proves
the claim in each degree.
\end{proof}

\begin{lemma}\label[lemma]{gl:reconstruction}
Let $n\geq1$, let $\pi$ be an irreducible ladder representation of $G_n$ on the
cuspidal line of $\chi$, and let $\Pi$ be a finite-length representation
of $G_n$ whose cuspidal support consists of characters of $G_1$.
Put $P=P_{1,n-1}$, with $r_P$ the identity when $n=1$.
If
\[
 (r_P\Pi)^{\mathrm{ss}}\simeq(r_P\pi)^{\mathrm{ss}},
\]
then $\Pi\simeq\pi$.
\end{lemma}
\begin{proof}
The case $n=1$ is immediate. By \eqref{eq:ladder-derivative}, the
constituents of $r_P\pi$ have pairwise distinct central characters.
Thus both Jacquet modules are semisimple and $r_P\Pi\simeq r_P\pi$.
We now use the
unique-socle property of \cite[Theorem~5.6]{Minguez2009}.
Choose an irreducible subrepresentation $\sigma\subseteq\Pi$.
Its Jacquet module is nonzero by the support hypothesis and is a
direct summand of $r_P\pi$. For an irreducible summand
$\rho\boxtimes\tau$, Frobenius reciprocity embeds both $\sigma$
and $\pi$ in $\rho\times\tau$. This induced representation has a
unique irreducible subrepresentation, so $\sigma\simeq\pi$.
Exactness gives $r_P(\Pi/\sigma)=0$. The support hypothesis then
forces $\Pi/\sigma=0$, proving the claim.
\end{proof}

\begin{proof}[Proof of \cref{thm:main}]
We first prove the equal-rank case $t=0$.
To apply \cref{gl:regular:theta}, list the lower endpoints in
increasing order as $(a_k,\ldots,a_1)$ and take $\sigma=w_0$.
The multisegment is balanced because its lower
endpoints decrease with its upper endpoints, excluding both forbidden
patterns in \cref{gl:regular:resolution-subsection}.
When the critical endpoints occur, their regular indices are
$i'=k+1-i$ and $j$, and the smaller permutation is $\sigma'=w'_0$.
Thus its degree is
\[
 \ell(w_0)-(2k-i'-j)-\ell(w'_0)=j-i=d_0(\fl),
\]
and its deleted multisegment is $\fl_0^*$.
For $d=1$, the nonsplit sequence \eqref{eq:regular-extension}
is exactly \eqref{eq:Sigma}. This proves all assertions for $t=0$.

The case $n=0$ follows from the rank-zero convention above.
Now assume $m>n\geq1$ and proceed by induction on $n+m$.
If the lower
endpoint $\nu^{s_t}$ or the upper endpoint $\nu^{-s_t}$ does not occur in $\fl$, we argue as in
\cite[Proposition~4.1]{CLLTZ2025}.

We thus may now assume \eqref{eq:unequal-critical}. Absorb $\chi$ into
the endpoints, so that $a_i=s_t$ and $b_j=-s_t$.
The ordinary recursion reduces $(n,m)$ to $(n-1,m-1)$, and the
exceptional recursion reduces it to $(n,m-1)$.
When $t=1$, the latter reaches the equal-rank case proved above.

Fix $r>0$. Let $\Pi^r_t(\fl)$ be $\Lrep(\fl_t^*)$ in the two cases
listed in \eqref{eq:unequal-positive}, and zero otherwise. We compare
its derivatives with those of $\cE^r_{n,m}(\Lrep(\fl))$.

First suppose $x\ne1-s_t$. By \eqref{eq:ladder-derivative}, a nonzero
$\partial_x\Lrep(\fl)$ is again a ladder representation. If
$x=b_j=-s_t$, lowering this upper endpoint removes the critical one,
and the corresponding positive-degree lift vanishes. The derivative of
$\Lrep(\fl_t^*)$ also vanishes, since $b_j$ has been deleted.
For the other removable upper endpoints, shortening commutes with
\eqref{eq:unequal-star}, and $d_t$ is unchanged whenever $d_t>0$.
To check the possible disappearance of a singleton, observe that
\[
 b_h-a_h+1\geq(-s_t+j-h)-(s_t-h+i)+1=d_t
 \qquad(i\leq h\leq j).
\]
Thus there is no such singleton if $d_t>1$. If $d_t=1$, a singleton
with $h<j$ forces $b_h=b_{h+1}+1$, so its upper endpoint cannot be shortened.
Singletons outside this interval do not change $j-i$. If $d_t\leq0$,
shortening an upper endpoint cannot increase $d_t$ as long as both critical
endpoints remain.

Conversely, the only new removable upper endpoint that deleting $b_j$ can create
is $b_{j-1}=1-s_t$, which is excluded here. Consequently, a standard
computation using \eqref{eq:ladder-derivative} and the induction
hypothesis gives
\[
 \partial_x \Pi^r_t(\fl)
 \simeq\cE^r_{n-1,m-1}(\partial_x\Lrep(\fl))
 \simeq\partial_x\cE^r_{n,m}(\Lrep(\fl)).
\]

It remains to consider $x=1-s_t$. By
\eqref{eq:exceptional-derivative-zero} and
\eqref{eq:exceptional-recursion},
\[
 \partial_{1-s_t}\cE^r_{n,m}(\Lrep(\fl))
 \simeq\nu^{1/2}\cE^r_{n,m-1}(\Lrep(\fl-1/2)),
\]
where $\fl-1/2=\sum_h[b_h-1/2,a_h-1/2]$. If there is no upper endpoint
$1-s_t$, the shifted ladder has no critical upper endpoint for the pair
$(n,m-1)$, and both sides of the desired comparison vanish. Otherwise
$b_{j-1}=1-s_t$. For the shifted ladder, the critical indices are
$i,j-1$, and
\[
 d_{t-1}(\fl-1/2)=(j-1)-i-(t-1)=d_t(\fl).
\]
Moreover, shortening the upper endpoint $b_{j-1}$ in \eqref{eq:unequal-star}
gives
\[
 \partial_{1-s_t}\Lrep(\fl_t^*)
 \simeq\nu^{1/2}\Lrep\bigl((\fl-1/2)_{t-1}^*\bigr).
\]
The induction hypothesis therefore proves the comparison for the
exceptional derivative as well.

Both $\cE^r_{n,m}(\Lrep(\fl))$ and $\Pi^r_t(\fl)$ have real-character
cuspidal support, by \cref{gl:support} and the definition of the candidate.
The preceding derivative equalities therefore give equal
semisimplifications of their $P_{1,m-1}$-Jacquet modules.
If $\Pi^r_t(\fl)\ne0$, it is a ladder and \cref{gl:reconstruction}
identifies the two representations. If $\Pi^r_t(\fl)=0$, the Jacquet
module of $\cE^r_{n,m}(\Lrep(\fl))$ vanishes, and its
character cuspidal support forces it to be zero. Taking duals gives
\eqref{eq:unequal-positive}.

Finally, the Euler--Poincar\'e formula in
\cref{gl:rank-homology}, followed by smooth duality, gives
\[
 \sum_{r\geq0}(-1)^r[\Theta_m^r(\Lrep(\fl))]
 =[\Lrep(\fl)^\vee\times\one_t].
\]
The positive-degree terms cancel when $d>1$, while for $d=1$ only
the term in degree one remains. This proves \eqref{eq:unequal-zero}
and the degree-zero assertion in the noncritical cases.
\end{proof}

\end{document}